\documentclass[11pt]{article}
\usepackage[a4paper, footskip=15mm,footnotesep=9mm,nohead,dvips]{geometry}
\usepackage{amsmath}
\usepackage{amsfonts}
\usepackage{amssymb}
\usepackage{amsthm}
\usepackage{epsfig, lscape}
\renewcommand{\leq}{\leqslant}
\renewcommand{\geq}{\geqslant}

\newtheorem{proposition}{Proposition}[section]
\newtheorem{corollary}{Corollary}[section]

\newtheorem{definition}{Definition}[section]
\newtheorem{rem}{Remark}[section]   
\newtheorem{exa}{Example}[section]   

\numberwithin{equation}{section}

\begin{document}

\title{A NOTE ON BRAID GROUP ACTIONS ON \\ SEMIORTHONORMAL BASES OF MUKAI LATTICES}
 \author{Amiel Ferman \footnote{Electronic mail: fermana1@math.biu.ac.il } \\
   Department of Mathematics, 
Bar-Ilan University \\ Ramat-Gan, Israel }
    \date{ }
 \maketitle

\begin{abstract}

We shed some light on the problem of determining the orbits of the braid group action on semiorthonormal bases of Mukai 
lattices as considered in \cite{GK04} and \cite{GO1}. We show that there is an algebraic (and in particular algorithmic) 
equivalence between this problem and the Hurwitz problem for integer matrix groups finitely generated by involutions. 
In particular we consider the case of $K_0(\mathbb P^n) \quad n \geq 4$ which was considered in \cite{GO1} and show
that the only obstruction for showing the transitivity of the braid group action on its semiorthonormal 
bases is the determination of the relations of particular finitely generated integer matrix groups. Although we prove transitivity
for an infinite set of Mukai lattices, our work, however, indicates quite strongly that the question of transitivity of semiorthonormal 
bases of Mukai lattices under the braid group action cannot be answered in general and can, at most, be resolved only in particular cases.

\end{abstract} 

%UDK: 514.85; MSC: 57R22


\section{Background and motivation} \label{sec: intro}

The bounded derived category of coherent sheaves of an algebraic variety $X$, $\mathcal D^b(X)$, has been
the focus of recent studies in various fields related with algebraic geometry (see for example \cite{B1}, \cite{BO1}).
Although $\mathcal D^b(X)$ in general contains less information than $Coh(X)$ (the abelian category of coherent sheaves on $X$),
it contains enough structure for the reconstruction of important invariants of $X$ such as higher Chow groups, K-theory and cohomology, as well
as birational geometric invariants of $X$. Furthermore, Bondal and Orlov \cite{BO1} have shown that projective varieties with ample 
canonical or anticanonical bundle are uniquely determined by their derived categories. 

The main idea of doing computations with $\mathcal D^b(X)$ consists of reducing the computations in the derived category essentially to linear algebra
using the Grothendieck $K_0$ functor. This was first conceived by Beilinson in \cite{BE}. We thus have a notion of 
a base, a so called \textbf{exceptional collection}, which consists of exceptional objects which generate $\mathcal D^b(X)$ and, by using $K_0$ becomes a 
vector space basis\footnote{strictly speaking it is a free $\mathbb Z$-module basis which becomes a vector space basis after tensoring with $\mathbb C$ for example} 
which is called a \textbf{semiorthonormal base}. For example, the exceptional collection of $K_0(\mathbb P^n)$ is induced
by $\{ \mathcal O, \mathcal O(1),\ldots,\mathcal O(n) \}$ (see \cite{BE,R1}). The exceptional collections of
del-Pezzo surfaces are also well studied \cite{KO}. We also have the \textbf{Euler form}, an (unimodular non-symmetric) inner product :

$$ \chi(E,F) = \sum (-1)^k dim Ext^k(E,F) $$ 
  
\leftline {defined on coherent sheaves and is respected by the $K_0$ functor.}
A semiorthonormal base equipped with such an inner product is called a \textbf{Mukai lattice}. 
Furthermore, there exists a non-trivial action of the braid group on such bases which is called \textbf{mutation} 
(we give the precise definitions in the next section). This action is transitive, for example, in the case of del-Pezzo surfaces \cite{KO}. 
Thus, in order to study semiorthonormal bases of Mukai lattices it is quite useful to study the orbits under the braid group action. In this paper
we address this problem which, as far as we know, was only considered in \cite{GK04} and in particular cases
in \cite{GO1} and \cite{R1}. For further details on this theory, which is known as \textbf{helix theory} one
can consult \cite{GK04} and \cite{R1}.

Our motivation in this paper was to try to understand the algebraic constraints in the constructions appearing in \cite{AKO1} and \cite{AKO2}.
In these works, Auroux Kuleshov and Orlov, have constructed categorical equivalences between bounded derived categories of certain 
algebraic varieties (or their deformations) and the so called \textbf{Fukaya-Seidel categories} of their corresponding dual symplectic manifolds according 
to the celebrated homological mirror symmetry conjecture of Kontsevich \cite{KON}. Apart from their physical significance (see for example \cite{HV}),
these constructions can also be viewed as a means for extracting effective algebraic invariants of symplectic manifolds. The Mukai lattices
can thus be considered as particular cases of \textbf{lagrangian vanishing cycles} appearing in Seidel's work \cite{SE1,SE2} on symplectic manifolds. Although these
constructions cannot be applied for general symplectic manifolds, they might give insights for the construction of analytic invariants
(as opposed to synthetic empirical invariants) of symplectic manifolds. We plan to pursue this point of view in future works. 

% \rtimes -- semidirect product

\section{Definitions and main results} 

We start off with the necessary definitions.

Recall that the {\bf braid group}, $B_n$, is a group with generators $\sigma_1,\ldots,\sigma_{n-1}$ and relations

$$ \sigma_i \sigma_{i+1} \sigma_i = \sigma_{i+1} \sigma_i \sigma_{i+1} \quad for \ i =1,\ldots,n-2 $$
$$ \sigma_i \sigma_j = \sigma_j \sigma_i \quad for \ |i-j| > 1 $$ 

\leftline{In what follows, $G$ is a group and $g_1,\ldots,g_n \in G$. }

\medskip

{\bf Hurwitz action} : This is an action of the braid group $B_n$ on $n$-tuples of group elements, defined on the 
generators of $B_n$ as follows (denote $g^h = hgh^{-1}$)

$$ \sigma_i (g_1,\ldots,g_n) = (g_1,\ldots,g_{i-1},g_{i+1}^{g_i}, g_i, g_{i + 2},\ldots,g_n) $$
$$ \sigma_i^{-1} (g_1,\ldots,g_n) = (g_1,\ldots,g_{i-1},g_{i + 1}, g_i^{g_{i+1}^{-1}}, g_{i + 2},\ldots,g_n) $$

\leftline{(that this is a well defined action see \cite{GK04}).}

\begin{rem} 

We shall sometimes call the Hurwitz action a \textbf{mutation} in accord with proposition \ref{EQ}
and the following definitions.
\end{rem}

$\medskip$

{\bf Hurwitz problem for groups} : Let $\{g_1,\ldots,g_n\}$ be an n-tuple of elements in the group $G$. Then the Hurwitz
problem for $G$ and $\{g_1,\ldots,g_n\}$  is to decide whether the tuple $(g_1,\ldots,g_n)$ and a given tuple $(h_1,\ldots,h_n)$
(where each $h_i$ is expressed in terms of the $g_j$-s or their inverses) are in the same
$B_n$-orbit, where the action of $B_n$ is the Hurwitz action defined above. 

\begin{rem}\label{AI} 

Note that our definition of the Hurwitz problem is more restrictive than \cite{TL05} where the
input is any two tuples of elements from the group. And so, the undecidability result, as well as the different
constructions in \cite{TL05} do not apply in our work.
\end{rem}

\bigskip

\begin{definition}\label{DEFMUK}
{\bf Mukai lattice} : This is a free $\mathbb Z$-module, $M$,  of finite rank, equipped with a bilinear map $M \times M \rightarrow \mathbb Z$
, denoted $\langle *,* \rangle$, which we assume to be {\bf unimodular}; this means that the natural map of $M$ into $Hom_{\mathbb Z}(M,\mathbb Z)$, sending
$m$ to $\langle m,* \rangle $ is a $\mathbb Z$-module isomorphism. A {\bf semi-orthonormal (or exceptional) basis} for a Mukai lattice
is a basis $\{E_1,\ldots,E_n\}$ of $M$ such that its Gram matrix $\chi_{ij} \stackrel{def}{=} \langle E_i, E_j \rangle$ satisfies
$\chi_{ii} = 1$ and $\chi_{ij} = 0$ for $i > j$ (We note that such a basis does not always exist but we shall work
with Mukai lattices which admit such bases). In this paper we shall use the term semi-orthonormal basis to mean
a tuple $(E_1,\ldots,E_n)$, where $\{E_1,\ldots,E_n\}$ is a semi-orthonormal basis of a Mukai lattice.
%Note that unimodularity means $det(\chi_{ij}) = \pm 1$ for all integer bases of $M$. 
Note that since $M$ is a 
free $\mathbb Z$-module, it can be faithfully represented as the free $\mathbb Z$-module $\mathbb Z^n$.
%In this paper (as done in \cite{GK04}) we shall consider the Mukai lattice to be a vector space over $\mathbb C$ after tensoring
%the module with $\mathbb C$, which clearly does not cause any loss of information.
\end{definition}

\bigskip

{\bf Mutations of semi-orthonormal bases} : This is an action of the braid group $B_n$ on semi-orthonormal bases
of a Mukai lattice of rank $n$ defined as follows

$$ \sigma_i (E_1,\ldots,E_n) = (E_1,\ldots,E_{i-1},L_{E_i}E_{i+1},E_i, E_{i+2},\ldots,E_n \} $$
$$ \sigma_i^{-1} (E_1,\ldots,E_n) = (E_1,\ldots,E_{i-1},E_{i+1}, R_{E_{i+1}}E_i, E_{i+2},\ldots,E_n \} $$

where 

$$ L_{E_i}E_{i+1} \stackrel{def}{=} E_{i+1} - \langle E_i,E_{i+1} \rangle E_i $$
$$ R_{E_{i+1}}E_i \stackrel{def}{=} E_i - \langle E_i,E_{i+1} \rangle E_{i+1} $$

\flushleft (that this is a well defined action on semi-orthonormal bases see \cite{GK04}).

\flushleft the actions of $\sigma_i$ and $\sigma_i^{-1}$ are called \bf left and right mutations \rm respectively. Thus, for example,
$L_{E_i}E_{i+1}$ is called the left mutation of $E_{i+1}$ by $E_i$.

\medskip

\leftline{Analogously to the Hurwitz problem on groups, we have}

\medskip

{\bf Hurwitz problem for Mukai lattices} : Given a Mukai lattice $M$ of rank $n$ and its semi-orthonormal
basis $(E_1,\ldots,E_n)$, the Hurwitz problem for $M$ is to decide whether a given semi-orthonormal
basis $(F_1,\ldots,F_n)$ is in the same $B_n$-orbit as $(E_1,\ldots,E_n)$ (for each $1 \leq i \leq n$, $F_i$ is given as a composition of left or right
mutations of $E_{\tau(i)}$ by the $E_i$-s where $\tau$ is a permutation of $1,\ldots, n$).

\medskip

As will be shown in the next proposition, the following set could be considered as equivalent, in some sense, to 
the set of all semi-orthonormal bases.

%\begin{definition}

%{\bf Reflect} : The set $Reflect_n$ contains, for every $n\times n$ symmetric integer matrix $B$ which has $2$ on its main diagonal, the tuple

%$$ (\varphi_1,\ldots,\varphi_n) $$

%where

%$$ \varphi_i(e_j) = e_j - B(i,j)e_i \quad \quad \varphi \in Hom_{\mathbb Z}(\mathbb Z^n) \quad e_i = (0,\ldots,0,1,0,\ldots,0) $$

%\flushleft In fact, each $\varphi_i$ is an automorphism of the $\mathbb Z$-module $\mathbb Z^n$ of order $2$ as can be readily seen
%by considering the $\varphi_i$ as lying in $Hom_{\mathbb C}(\mathbb C^n)$ for example. Note also that the $\varphi_i$ in 
%each tuple generate a subgroup of $Aut_{\mathbb Z}(\mathbb Z^n)$.
%\end{definition}

\begin{definition}\label{DEFREF}

{\bf Reflect} : Given a Mukai lattice $M$, 
the set $Reflect(M)$ consists of the tuples $(\psi_{E_1},\ldots,\psi_{E_n})$ of elements in $Hom_{\mathbb Z}(M)$ 
where $\{E_1,\ldots, E_n\}$ is a semi-orthonormal base of $M$ and

$$ \psi_{E_i}(E_j) = E_j - (\langle E_i, E_j \rangle + \langle E_j, E_i \rangle)\cdot E_i $$ 

\flushleft we shall see below that in fact each $\psi_{E_i}$ is an automorphism of $M$ of order $2$.

\end{definition}

\medskip

\begin{definition}
{\bf M-Bases} Given a Mukai lattice $M$, the set $M-Bases(M)$ consists of the tuples $(E_1,\ldots,E_n)$
where $\{E_1,\ldots,E_n\}$ is a semi-orthonormal base of $M$.

\end{definition}

Mutations of semi-orthonormal bases of a Mukai lattice can be considered in terms of Hurwitz action on group tuples as shown in the following

\begin{proposition}[Bijective Equivariance]\label{EQ}

Given a Mukai lattice $M$ of rank $n$, there exists a bijective equivariant mapping between the set $M-Bases(M)$ and the set 
$Reflect(M)$ w.r.t the $B_n$ action up to the sign of the vectors in the semi-orthonormal bases.

\end{proposition}
\begin{proof}

\bigskip

Given a Mukai lattice $M$ of rank $n$ with a semiorthonormal base $(E_1,\ldots,E_n)$ and an inner product $\langle \cdot,\cdot \rangle$,
first let us denote its symmetrized inner product by $B$, that is

$$ B(E,F) \stackrel{def}{\equiv} \langle E,F \rangle + \langle F,E \rangle $$

Consider the following map
 from $M-Bases(M)$ into $Reflect(M)$

$$ \Psi (E_1,\ldots,E_n) \stackrel{def}{=} (\psi_{E_1},\ldots,\psi_{E_n}) $$

where we define $\psi_{E_i} \in Hom_{\mathbb Z}(M)$ as follows

$$ \psi_{E_i}(F) = F - B(E_i,F)E_i \quad F,E_i \in M $$

First, let us show that $\psi_{E_i}$ is a map of order $2$ for all $i=1\ldots n$ :

$$ \psi_{E_i}(\psi_{E_i}(A)) = \psi_{E_i}(A - B(E_i,A)E_i) = A - B(E_i,A)E_i - B(A - B(E_i,A)E_i,E_i)E_i = $$
$$ = A - B(E_i,A)E_i - B(A,E_i)E_i + B(E_i,A)\cdot B(E_i,E_i)E_i = A $$

%\leftline{(Recall that $B(E_i,E_i) = \langle E_i, E_i \rangle + \langle E_i, E_i \rangle = 2$)}

\bigskip

Secondly, let us show that $\psi_{E_i}$ is an isometry w.r.t $B$ :

\begin{gather}
B(\psi_{E_i}(F),\psi_{E_i}(G)) = B(F - B(F,E_i)E_i, G - B(G,E_i)E_i) = \notag \\
= B(F,G) - B(G,E_i)\cdot B(F,E_i) - B(F,E_i)\cdot B(E_i,G) + B(F,E_i)\cdot B(G,E_i) \cdot B(E_i,E_i) = \notag \\ 
= B(F,G) - 2B(G,E_i)\cdot B(F,E_i) + 2B(F,E_i)\cdot B(G,E_i) = \notag \\ 
= B(F,G)\label{ISO}
\end{gather}

\bigskip

Now, that $\Psi$ is surjective is clear from definition \ref{DEFREF}. To show it is injective up to the signs of the $E_i$-s note first that

$$ \psi_{-E_i}(F) = F - B(-E_i,F)\cdot(-E_i) =  F - B(E_i,F)E_i = \psi_{E_i}(F) $$

\flushleft on the other hand, suppose that $(E_1,\ldots,E_n)$ and $(F_1,\ldots,F_n)$ are two semi-orthonormal bases of $M$ and
suppose that

$$ (\psi_{E_1},\ldots,\psi_{E_n}) = (\psi_{F_1},\ldots,\psi_{F_n}) $$

\flushleft then for every $i$ we have that 

$$ \psi_{E_i}(E_i) = \psi_{F_i}(E_i) \Longrightarrow -E_i = E_i - B(F_i,E_i)F_i \Longrightarrow 2E_i =   B(F_i,E_i)F_i $$

$$\psi_{E_i}(F_i) = \psi_{F_i}(F_i) \Longrightarrow F_i - B(F_i,E_i)E_i = -F_i \Longrightarrow 2F_i =  B(F_i,E_i)E_i $$

\flushleft the last two equations imply that $4E_i = (B(E_i,F_i))^2 E_i$ but this must mean that $B(E_i,F_i) = \pm 2$ 
and hence $E_i = \pm F_i$.

\medskip

%Thus, the hyperplane $span\{F_j^i \}_{j\neq i}$ is fixed by $\psi_{E_i}$, and $F_i^i = E_i$ 
%is inverted by $\psi_{E_i}$. This proves that $\psi_{E_i}$ is an isometry (w.r.t $B$) of order 2. We shall therefore write from now on $\psi_{E_i}$ instead of $\psi_{E_i}^{-1}$.

%First, let us show that $\psi_{E_i}$ is a map of order $2$ for all $i=1\ldots n$ :

%$$ \psi_{E_i}(\psi_{E_i}(A)) = \psi_{E_i}(A - B(E_i,A)E_i) = A - B(E_i,A)E_i - B(A - B(E_i,A)E_i,E_i)E_i = $$
%$$ = A - B(E_i,A)E_i - B(A,E_i)E_i + B(E_i,A)\cdot B(E_i,E_i)E_i = A $$

%\leftline{(Recall that $B(E_i,E_i) = \langle E_i, E_i \rangle + \langle E_i, E_i \rangle = 2$)}

\bigskip

%Secondly, let us show that $\psi_{E_i}$ is an isometry w.r.t $B$ :

%\begin{gather}
%B(\psi_{E_i}(F),\psi_{E_i}(G)) = B(F - B(F,E_i)E_i, G - B(G,E_i)E_i) = \notag \\
%= B(F,G) - B(G,E_i)\cdot B(F,E_i) - B(F,E_i)\cdot B(E_i,G) + B(F,E_i)\cdot B(G,E_i) \cdot B(E_i,E_i) = \notag \\ 
%= B(F,G) - 2B(G,E_i)\cdot B(F,E_i) + 2B(F,E_i)\cdot B(G,E_i) = \notag \\ 
%= B(F,G)\label{ISO}
%\end{gather}

%\bigskip

%Now, let us show that $\Psi$ is well defined. We have that $\psi_{E_i}(E_i) = -E_i $ and for all $j\neq i$, denote $F_j^i = \psi_{E_i}E_j + E_j$, then

%So we now have that, for all $i$, $\psi_{E_i} \in$ Reflect($M$) and $\Psi$ is indeed well defined.

\bigskip

In what follows, we shall write $\psi_{E_i}$ instead of $\psi_{E_i}^{-1}$, since, as mentioned above $\psi_{E_i}$ is a map of order $2$.

\bigskip

To show equivariance, let us first show that for each $i$, $1 \leq i < n$, we have

\begin{equation}
\psi_{L_{E_i}E_{i+1}} = \psi_{E_{i+1}}^{\psi_{E_i}} ( = \psi_{E_i} \psi_{E_{i+1}} \psi_{E_i} ) 
\label{EXCON}
\end{equation}

To prove this, we shall show that both sides of equation \ref{EXCON} agree on a basis of $M$. Consider then the following
set of vectors $\{ F_k^j \}_{k=1\ldots n}$ which is defined for each $j$, $1 \leq j \leq n$, as follows :

\begin{equation}
 F_k^j = \begin{cases}
\psi_{E_j}E_k + E_k & k\neq j \\
E_j & k = j
\end{cases} 
\label{FDEF}
\end{equation}

\flushleft that this is indeed a basis for $M$ is easy to see by considering the faithful representation of $M$ in $\mathbb Z^n$
as mentioned in definition \ref{DEFMUK} where the $E_i$-s are mapped to the standard basis vectors : After this representation is applied it is 
easy to see that the matrix whose columns are $F^j_k$ has determinant $1$ and hence it is inverted by an integer matrix which
means that its columns span $\mathbb Z^n$ over $\mathbb Z$.  Now, for a fixed $i$, we have that $\{\psi_{E_i}(F_k^{i+1}) \}_{k=1\ldots n}$
is also a basis of $M$ (as each $\psi_{E_i}$ is an automorphism of $M$). We shall thus show that both sides of equation \ref{EXCON}
agree on the basis $\{\psi_{E_i}(F_k^{i+1}) \}_{k=1\ldots n}$ of $M$.

%the set of vectors $\{ \psi_{E_i}(F_k^j) \}_{k = 1 \ldots n}$ 
%which is a basis since it is the application of the isometry $\psi_{E_i}$ on the basis $\{ F_k^j \}_{k = 1\ldots n} $.  

\bigskip

Note first that

\begin{equation}
\forall A \in M, \ \forall i=1\ldots n : \quad B(A,E_i) = 0 \ \Longleftrightarrow \ \psi_{E_i}(A) = A
\label{B0}
\end{equation}

and that for all $k$, $k\neq i$ :

\begin{equation}
\psi_{E_i}(F_k^i) = 
\psi_{E_i}(\psi_{E_i}E_k + E_k) = 
\psi_{E_i}(\psi_{E_i}(E_k)) + \psi_{E_i}(E_k) =   E_k + \psi_{E_i}(E_k) = F_k^i 
\label{FIX}
\end{equation}

\bigskip

Now, on the one hand, using equation \ref{FIX}, we have that 

$$ \psi_{E_{i+1}}^{\psi_{E_i}} (\psi_{E_i}(F_k^{i+1})) =  \psi_{E_i} \psi_{E_{i+1}} \psi_{E_i} (\psi_{E_i}(F_k^{i+1})) = \psi_{E_i}(F_k^{i+1}) \quad k\neq {i+1} $$

$$ \psi_{E_{i+1}}^{\psi_{E_i}} (\psi_{E_i}(F_{i+1}^{i+1})) = \psi_{E_{i+1}}^{\psi_{E_i}} (\psi_{E_i}(E_{i+1})) =  \psi_{E_i} \psi_{E_{i+1}} \psi_{E_i} (\psi_{E_i}(E_{i+1})) = -\psi_{E_i}(E_{i+1}) = -\psi_{E_i}(F_{i+1}^{i+1}) $$

\leftline{and on the other hand we have that}

$$ B(L_{E_i}E_{i+1}, \psi_{E_i}(F_k^{i+1})) = B(\psi_{E_i}(E_{i+1}), \psi_{E_i}(F_k^{i+1})) = B(E_{i+1}, F_k^{i+1}) = 0 $$

the second equality follows the fact that $\psi_{E_i}$ is an isometry and the third equality follows \ref{FIX} and \ref{B0}; now, using \ref{B0} we have that this equality is equivalent to

$$ \psi_{L_{E_i}E_{i+1}}(\psi_{E_i}(F_k^{i+1})) = \psi_{E_i}(F_k^{i+1}) \quad k\neq {i+1} $$

\flushleft and furthermore we have

$$ \psi_{L_{E_i}E_{i+1}}(\psi_{E_i}(F_{i+1}^{i+1})) = \psi_{L_{E_i}E_{i+1}}(\psi_{E_i}(E_{i+1})) = \psi_{\psi_{E_i}(E_{i+1})}(\psi_{E_i}(E_{i+1})) = -\psi_{E_i}(E_{i+1}) = -\psi_{E_i}(F_{i+1}^{i+1})         $$

Thus, we have shown that both sides of equation \ref{EXCON} agree on the basis $\{\psi_{E_i}(F_k^{i+1}) \}_{k=1\ldots n}$ , and so

\begin{equation}
\psi_{L_{E_i}E_{i+1}} = \psi_{E_{i+1}}^{\psi_{E_i}}
\label{LE}
\end{equation}

\leftline{as $R_{E_{i+1}}(E_i) = \psi_{E_{i+1}}(E_i)$, a completely analogous computation shows that}

\begin{equation}
\psi_{R_{E_{i+1}}E_i} = \psi_{E_{i+1}}^{\psi_{E_i}}
\label{RE}
\end{equation}

Now, \ref{LE} and \ref{RE} show that for every braid group generator $\sigma_i$ we have

$$ \Psi(\sigma_i(E_1,\ldots,E_n)) = \sigma_i(\Psi(E_1,\ldots,E_n)) $$

\flushleft proceeding with induction on the number of generators this establishes equivariance.

\end{proof}

\bigskip

\begin{rem}\label{CONJ}
There exists another natural action on semi-orthonormal bases. Consider $Isom(M)$ : The group of isometries of $M$ as an inner
product space. Then the action of $\varphi \in Isom(M)$ on the semiorthonormal 
base $(E_1,\ldots,E_n)$ is $(\varphi(E_1),\ldots,\varphi(E_n))$. As can be seen by the following computation :

$$ \psi_{\varphi(E_i)}(\varphi(E_j)) = \varphi(E_j) - B(\varphi(E_i),\varphi(E_j))\varphi(E_i) = \varphi\circ \psi_{E_i}(E_j) $$

and hence 

$$ \psi_{\varphi(E_i)}\circ \varphi = \varphi \circ \psi_{E_i} $$

hence

$$ \psi_{E_i}^{\varphi} = \psi_{\varphi(E_i)} \quad \varphi \in Isom(M) $$

\flushleft and so, in the notation of proposition \ref{EQ}, we have that 

$$ \Psi(\varphi\cdot(E_1,\ldots,E_n)) = (\psi_{E_1}^\varphi,\ldots,\psi_{E_n}^\varphi) \quad \varphi \in Isom(M) $$

\flushleft  so the action of $Isom(M)$ on the semiorthonormal bases corresponds to \textbf{global conjugtion}
on the corresponding group tuples. We will address this subject more generally in the last section.

(In \cite{GK04} the considered action on semiorthonormal bases of Mukai lattices is actually the action of the group $Isom(M) \rtimes B_n$,
however, it is not clear what is meant by this semi-product as the stabilizer of $B_n$ for example may be non-trivial. For this reason we decided
to address the action of $Isom(M)$ separately).

\end{rem}

Proposition \ref{EQ} gives a good notion of the algebraic difficulty of analysing the orbits
in Mukai lattices. Indeed, although the groups generated by the tuples in $Reflect(M)$ are finitely generated, there is no reason
that they should be finitely presented (i.e. have a finite number of relations). However, as the following
proposition shows, groups which have simple relations (such as free Coxeter groups) are easy to analyse 
in this context. 

Before stating the proposition, let us first say what we mean by a \textbf{transitive} orbit of the $B_n$ action
on a group tuple $(g_1,\ldots,g_n)$. Note that if $(h_1,\ldots, h_n)$ is in 
the same $B_n$ orbit of $(g_1,\ldots, g_n)$ (where each $h_i$ is expressed in terms of the $g_i$) then necessarily $h_i = g_{\tau(i)}^{f_i}$ where $\tau$ is a permutation 
of $(1,\ldots,n)$ and $f_i \in G$. Furthermore, we have $h_1\cdots h_n = g_1 \cdots g_n$. 
Both properties are preserved under a single mutation, as can be readily checked, and in general follow by induction on the number of mutations.
We shall thus say that the $B_n$ action on $(g_1,\ldots,g_n)$ is \textbf{transitive} if every tuple of the form $(g_{\tau(1)}^{f_1},\ldots,g_{\tau(n)}^{f_n})$
is in the orbit of $(g_1,\ldots,g_n)$.

\begin{proposition}\label{FREE}
Let $G = \{ g_1,\ldots,g_n | g_i^2 = 1, \quad i = 1\ldots n \}$ (i.e. $G$ is a free Coxeter group),
then the $B_n$ action on $(g_1,\ldots,g_n)$ is transitive.
\end{proposition}
\begin{proof}

Throughout the proof we shall assume that all elements of $G$ are presented in terms of the generating set $\{g_1,\ldots,g_n\}$.

To prove the claim, we shall show that, given a tuple $(g_{\tau(1)}^{f_1},\ldots,g_{\tau(n)}^{f_n})$, 
where $\tau$ is a permutation of $(1,\ldots,n)$ and $f_i \in G$ and furthermore  $g_{\tau(1)}^{f_1} \cdots g_{\tau(n)}^{f_n} = g_1\cdots g_n$, we have that
$(g_{\tau(1)}^{f_1},\ldots,g_{\tau(n)}^{f_n})$ is in the same $B_n$ orbit of $(g_1,\ldots,g_n)$. 

\bigskip

Since $G$ is a free Coxeter group, we can use the notion of length of words and its properties (see \cite{H90} chapter 5).
We denote the length of a word $u$, as usual, by $l(u)$.

\bigskip 

We shall prove that given
a tuple $(g_{\tau(1)}^{f_1},\ldots,g_{\tau(n)}^{f_n})$ there always exists a mutation (at least one such) that reduces
the total length of words in $(g_{\tau(1)}^{f_1},\ldots,g_{\tau(n)}^{f_n})$  . This shall be proved by induction on the total length of 
the words in $(g_{\tau(1)}^{f_1},\ldots,g_{\tau(n)}^{f_n})$. We assume that each word in the given tuple is reduced.

\bigskip

The base of the induction is simply the tuple $(g_1,\ldots,g_n)$ and the claim follows trivially.

%induction base
%$$ (g_{\tau(1)}^{f_1},\ldots,g_{\tau(n)}^{f_n}) = (g_{\tau(1)}^{g_j},g_{\tau(2)},\ldots,g_{\tau(n)}) $$

%then we have that

%$$ g_j g_{\tau(1)}g_j  g_{\tau(2)} \cdots g_{\tau(n)} = g_1\cdots g_n  $$

%in a free Coxeter group this equation implies

%$$ g_i = g_{\tau(i)} \quad i=3\ldots n $$

%and

%$$ g_j g_{\tau(1)}g_j g_{\tau(2)} g_2 g_1 = 1  $$ 

%now, if $\tau(2) = 2$ then $\tau(1) = 1$ and $j=1$ which contradicts our assumption 
%that each word is reduced, so we must have that $\tau = (1 \ 2)$ and $j=1$ and so

%$$ (g_{\tau(1)}^{g_j},g_{\tau(2)},\ldots,g_{\tau(n)}) = (g_2^{g_1},g_1,g_3,\ldots,g_n) = \sigma_1(g_1,\ldots,g_n) $$

%which proves the induction base.

To prove the induction step, we will show that given the tuple

$$ (g_{\tau(1)}^{f_1},\ldots,g_{\tau(n)}^{f_n}) $$

\flushleft there exists a mutation for this tuple which reduces the total length of the words (and hence the claim will follow by the induction hypothesis).

\bigskip

In the following equations we shall use the fact that for every word $u$ in a Coxeter group $l(u) = l(u^{-1})$. Furthermore, we
shall use the following computation (recall that we denote $g^h = hgh^{-1}$)

\begin{equation}
(g^h)^f = fhgh^{-1}f^{-1} = g^{fh}
\label{SHORT}
\end{equation}

Now, applying a left mutation on the given tuple $ (g_{\tau(1)}^{f_1},\ldots,g_{\tau(n)}^{f_n}) $ gives

\begin{equation}
(g_{\tau(1)}^{f_1},\ldots,  g_{\tau(i-1)}^{f_{i-1}},  g_{\tau(i+1)}^{f_{i+1}}, (g_{\tau(i)}^{f_i})^{( g_{\tau(i+1)}^{f_{i+1}} )^{-1}} , g_{\tau(i+2)}^{f_{i+2}},\ldots, g_{\tau(n)}^{f_n}) 
\label{LEFTM}
\end{equation}

And, similarly, applying a right mutation on the given tuple $ (g_{\tau(1)}^{f_1},\ldots,g_{\tau(n)}^{f_n}) $ gives

\begin{equation}
(g_{\tau(1)}^{f_1},\ldots,  g_{\tau(i-1)}^{f_{i-1}}, (g_{\tau(i+1)}^{f_i+1})^{g_{\tau(i)}^{f_{i}}}, g_{\tau(i)}^{f_i}, g_{\tau(i+2)}^{f_{i+2}},\ldots, g_{\tau(n)}^{f_n})
\label{RIGHTM}
\end{equation}

Note that the only new element after applying a left mutation in the $i$-th place on the given tuple is

\begin{equation} 
(g_{\tau(i)}^{f_i})^{( g_{\tau(i+1)}^{f_{i+1}} )^{-1}} = (g_{\tau(i)}^{f_i})^{g_{\tau(i+1)}^{f_{i+1}}}  
=  (g_{\tau(i)})^{g_{\tau(i+1)}^{f_{i+1}}\cdot f_i}  =  (g_{\tau(i)})^{f_{i+1} g_{\tau(i+1)} f_{i+1}^{-1}  f_i}
\label{NEWR}
\end{equation}

Similarly, the only new element after applying a right mutation in the $i$-th place is

\begin{equation}
(g_{\tau(i+1)}^{f_i+1})^{g_{\tau(i)}^{f_{i}}} =  (g_{\tau(i+1)})^{g_{\tau(i)}^{f_{i}}\cdot f_{i+1}}  
=  (g_{\tau(i+1)})^{f_{i} g_{\tau(i)} f_{i}^{-1}  f_{i+1}}
\label{NEWL}
\end{equation}

The change of the total length of the words comprising the tuple in case of the left mutation is

\begin{equation}
l((g_{\tau(i)})^{g_{\tau(i+1)}^{f_{i+1}}\cdot f_i}  ) - l((g_{\tau(i)}^{f_i}))
\label{CHR}
\end{equation}

and in the case of a right mutation it is 

\begin{equation}
l((g_{\tau(i+1)})^{g_{\tau(i)}^{f_{i}}\cdot f_{i+1}}  ) - l((g_{\tau(i+1)}^{f_i+1}))
\label{CHL}
\end{equation}

We would like to show that either of these quantities is negative for some $i$, hence it is enough to show that there exists at least one $i$ for which either

\begin{equation}
l(f_{i}^{-1} g_{\tau(i+1)}^{f_{i+1}}) < l(f_i)
\label{RESR}
\end{equation}

or

\begin{equation}
l(g_{\tau(i)}^{f_{i}} f_{i+1}) < l(f_{i+1})
\label{RESL}
\end{equation}

Consider the following equation which holds by assumption

$$ g_{\tau(1)}^{f_1} \cdots g_{\tau(n)}^{f_n} = g_1\cdots g_n $$

or more explicitely

$$ f_1 g_{\tau(1)} f_{1}^{-1} f_2 g_{\tau(2)} f_{2}^{-1}  \cdots f_n g_{\tau(n)} f_{n}^{-1}= g_1\cdots g_n   $$

If none of the $g_{\tau(i)}$ appearing on the left hand side of this equality is cancelled in order for the
equation to hold then necessarily for all $i$ we have that $l(f_i) = 0$ (i.e. we are in the induction base). 
So we can assume that for some $i$,  $g_{\tau(i)}$ is cancelled during the cancellation process. We can assume,
without lose of generality, that the subword which contains the letter that cancels $g_{\tau(i)}$ is on its right. 
Consider then the following subword appearing on the left hand side of the above equation 

$$ u = g_{\tau(i)} \cdot f_{i}^{-1} \cdot f_{i+1} \cdots f_{k-1}^{-1} \cdot f_k \cdot g_{\tau(k)} \cdot f_{k}^{-1} \cdot f_{k+1} \cdots f_{j-1}^{-1} \cdot f_j \cdot g_{\tau(j)} f_j^{-1}$$

\flushleft in this subword we choose $k$ to be maximal such that $v = g_{\tau(i)} \cdot f_{i}^{-1} \cdot f_{i+1} \cdots f_{k-1}^{-1}$ is reduced.
Note that the last letter of the reduced subword $v$ in $u$ must be the last letter of $f_l^{-1}$ for some $l$ as all the subwords
of the form $f_l\cdot g_{\tau(l)} \cdot f_l^{-1}$ are reduced by assumption. Furthermore, we choose $j$ to be the minimal index such that 
$g_{\tau(j)}$ does not appear in $v^{-1}$, or we choose $j=n$ if all $g_{\tau(j)}$ appear in $v^{-1}$ for all $j > i$. 
In any case, as $\tau$ is a permutation of $(1 \ldots n)$, we must have that the last word of $v^{-1}$, which is $g_{\tau(i)}$,
appears in $f_{j-1}^{-1}$ or in $f_j$ or in $f_n^{-1}$. Consider the following subword in $u$ (for the same index $k$ as above)

\begin{equation}
g_{\tau(k-1)} \cdot f_{k-1}^{-1} \cdot f_k \cdot g_{\tau(k)}
\label{EXPK}
\end{equation}

\flushleft clearly, either of $f_{k-1}^{-1}$ or $f_k$ is not trivial (for otherwise $g_{\tau(k-1)}$ and $g_{\tau(k)}$ could not 
be cancelled). Now, first we address the case where both $ g_{\tau(k-1)}$ and $g_{\tau(k)} $ are cancelled. Suppose first that $l(f_k) > l(f_{k-1}^{-1})$ ($l(f_k) \neq l(f_{k-1}^{-1})$ since
$g_{\tau(k-1)}$ cannot cancel $g_{\tau(k)}$). Then necessarily $f_k = (g_{\tau(k-1)} \cdot f_{k-1}^{-1})^{-1} \cdot w$  for some reduced (possibly trivial) word $w$. 
For otherwise $g_{\tau(k-1)}$ would remain uncancelled as opposed to the assumption. Similarly, if $l(f_k) < l(f_{k-1}^{-1})$
then  $f_{k-1}^{-1} = y\cdot (f_k \cdot g_{\tau(k)})^{-1} $ for some reduced (possibly trivial) word $y$.  In the first case we have that

$$ l(f_{k-1} g_{\tau(k-1)} f_{k-1}^{-1} f_k) \leq l(f_k) - l( f_{k-1}^{-1}) - l( g_{\tau(k-1)}) + l(f_{k-1}) = $$ 
$$ = l(f_k) - l( f_{k-1}) - 1 + l(f_{k-1}) < l(f_k)    $$

Thus we have that inequality \ref{RESL} holds for $i=k-1$. Similarly, in the second case we have that inequality \ref{RESR} holds for $i=k$.

\bigskip

Now, while $g_{\tau(k-1)}$ allways cancels in the expression \ref{EXPK} according to our choice of indices, $g_{\tau(k)}$ may not cancel,
however, we must have in this case that $f_k = v^{-1}\cdot w$ for some reduced (possibly trivial) word $w$. Hence this must be 
the first case considered in the above argument (i.e. that  $f_k = (g_{\tau(k-1)} \cdot f_{k-1}^{-1})^{-1} \cdot w$  for some reduced (possibly trivial) word $w$)
and so the claim follows.

\end{proof}

\begin{exa}\label{K0P4}
In light of proposition \ref{FREE}, let us consider the case of $K_0(\mathbb P^4)$ which
is yet to be determined (see \cite{GO1}). Its Gram matrix is easily seen to be

$$ \left( \begin{matrix}
1 & 5 & 15 & 35 & 70\\
0 & 1 & 5 & 15 & 35 \\
0 & 0 & 1 & 5 & 15 \\
0 & 0 & 0 & 1 & 5 \\
0 & 0 & 0 & 0 & 1
\end{matrix}\right) $$

\flushleft Let us write $Reflect(K_0(\mathbb P^4))$ explicitly. It is generated by the following five matrices :

$$ s_1 = \left( \begin{matrix}
-1 & -5 & -15 & -35 & -70\\
0 & 1 & 0 & 0 & 0 \\
0 & 0 & 1 & 0 & 0 \\
0 & 0 & 0 & 1 & 0 \\
0 & 0 & 0 & 0 & 1
\end{matrix}\right) \quad
s_2 = \left( \begin{matrix}
1 & 0 & 0 & 0 & 0 \\
-5 & -1 & -5 & -15 & -35\\
0 & 0 & 1 & 0 & 0 \\
0 & 0 & 0 & 1 & 0 \\
0 & 0 & 0 & 0 & 1
\end{matrix}\right) $$

$$ s_3 = \left( \begin{matrix}
1 & 0 & 0 & 0 & 0 \\
0 & 1 & 0 & 0 & 0 \\
-15 & -5 & -1 & -5 & -15\\
0 & 0 & 0 & 1 & 0 \\
0 & 0 & 0 & 0 & 1
\end{matrix}\right) \quad
s_4 = \left( \begin{matrix}
1 & 0 & 0 & 0 & 0 \\
0 & 1 & 0 & 0 & 0 \\
0 & 0 & 1 & 0 & 0 \\
-35 & -15 & -5 & -1 & -5 \\
0 & 0 & 0 & 0 & 1
\end{matrix}\right) $$

$$ s_5 = \left( \begin{matrix}
1 & 0 & 0 & 0 & 0 \\
0 & 1 & 0 & 0 & 0 \\
0 & 0 & 1 & 0 & 0 \\
0 & 0 & 0 & 1 & 0 \\
-70 & -35 & -15 & -5 & -1
\end{matrix}\right) $$

According to proposition \ref{FREE}, the action of the braid group on semiorhonormal bases of $K_0(\mathbb P^4)$
is transitive if the group generated by $\{s_1,s_2,s_3,s_4,s_5 \}$ is a free Coxeter group. It is not hard to check that
this group contains $\mathbb F_2$ (the free group with $2$ generators) according to the Tits alternative \cite{B99}, however
it seems quite difficult to determine the relations in this group\footnote{Stefan Kohl who helped with the analysis of this group
determined that there are strong 'empricial' evidence that this group is indeed a free Coxeter group, but a proof (for example
using the table-tennis lemma) seems hard}. 
\end{exa}

$$ $$

Finally let us shortly discuss the algorithmic aspect of the problem. Recall that an algorithmic problem $A_1$ is said to be (Turing-) Reducible to an algorithmic problem $A_2$
if, given an algorithm that solves $A_2$, the problem $A_1$ can also be solved. Furthermore $A_1$ is (Turing-)equivalent to $A_2$ 
if each of them is reducible to the other. Hence, proposition \ref{EQ} immidately gives

\begin{corollary}\label{ALEQ}
For a given Mukai lattice $M$, the Hurwitz problem for Reflect($M$) is equivalent to the Hurwitz problem for $M$.
\end{corollary}

This observation is useful especially on the negative side : it is sometimes easier to see that an algorithmic
problem is undecidable rather than analysing it algebraically. In our specific case it might be easier
to determine such undecidability results for the Hurwitz problem on $Reflect(M)$ rather than trying to disprove the
transitivity of the action on the Mukai lattice.

$$ $$

\section{Conclusions and further remarks}

The bijective equivariance between Mukai lattices and group tuples of reflections proved in proposition \ref{EQ} indicates that the main 
difficulty in analysing the braid group orbits lies in the analysis of the relations in those groups of reflections. These
groups are matrix integer groups finitely generated by matrices of order $2$. It seems that there exists no theory that could determine
in general the relations in these kind of groups. Indeed, as noted above, there is no guarntee that these groups are 
even finitely presented (i.e. have a finite number of relations). As can be seen in example \ref{K0P4}, even particular cases
seem quite hard to analyse. Although we do not have a specific example, it is quite reasonable to assume
that there exists a particular Mukai lattice for which the Hurwitz problem is undecidable (see for example the diagram 
on page 31 in \cite{M2} for the different undecidability results of finitely generated linear groups). Furthermore,
it would also be quite reasonable to assume that the question of transitivity is not algorithmic in general.

According to remark \ref{CONJ}, the conjugacy problem for a group of the form $Reflect(M)$ is algorithmically reducible
to the problem of determining whether two semiorthonormal bases are in the same $Isom(M)$ orbit. Thus, the undecidability
of the conjugacy problem for $Reflect(M)$ would imply the undecdability for the $Isom(M)$ orbit problem. This fact,
combined with the different undecidability results for finitely generated linear groups in \cite{M1} and \cite{M2}
might imply that it would be even easier to construct such an undecidability result in this case.

A natural extension of proposition \ref{FREE} (apart from the obvious extension to general Coxeter groups) might be
to prove the same claim for \textbf{biautomatic groups} (a free Coxeter group is an example of a biautomatic group),
intuitively these groups have the property that relations in them have a computationaly simple description which
allows for the solution of different decision problems defined on them (see section 7 in \cite{M2}). However, as noted
above, it is not reasonable to assume that these kind of results could be streangthened to provide a sufficient condition
for the transitivity of the braid action.

Finally, although our work implies that the analysis of braid group actions on (semiorthonormal bases of) general Mukai lattices is not algorithmic
in general, it is not clear wheather the same holds for Mukai lattices induced from algebraic varieties.
It would be interesting to understand exactly which Mukai lattices can or cannot be induced from algebraic varieties.
Although this problem also seems quite difficult (if not altogether impossible in general), such an understanding could imply 
at least which kind of Mukai lattices are worth focusing our attention on. 

\bigskip 

{\bf Acknowledgements}. The author would like to thank : Vadim Ostapenko, Yann Sepulcre and Maxim Leyenson for helpful discussions
in algebraic geometry and in general. Bettina Eick, Stefan Kohl and the group-pub-forum for their help with different questions in group theory.
And of course my two advisors : Tahl Nowik and Mina Teicher for their help and support.

\end{document}